\documentclass{article}

\usepackage{amsmath}
\usepackage{url}

\begin{document}

\title{A Proof That Zeilberger Missed:  A New Proof Of An Identity By Chaundy
And Bullard Based On The Wilf-Zeilberger Method}
\author{YiJun Chen}
\date{}
\maketitle

\begin{abstract}
In this paper, based on the Wilf-Zeilberger theory, a very succinct new proof of
an identity by Chaundy and Bullard is given.
\end{abstract}

In their paper \cite{koornwinder}, T. H. Koornwinder and M. J. Schlosser have given detailed
explanations of the historical developments, different proofs (in fact, they provide seven different proofs), and different generalizations of the identity

\begin{equation}
\label{eq1}
1 = (1 - x)^{n + 1}\sum\limits_{k = 0}^m {{n+k} \choose {k}} x^{k} + x^{m + 1}\sum\limits_{k = 0}^n {{m+k} \choose {k}}(1 - x)^{k},
\end{equation}

\noindent
which is attributed to Chaundy and Bullard. In particular, we can learn in
\cite{koornwinder} that the special case where $m = n$ plays a key role in Daubechies's
theory of wavelets of compact support(see \cite{daubechies1}or \cite{daubechies2}). A very succinct
probabilistic proof of this special case was given by Zeilberger (see \cite{zeilberger}).
As pointed out in \cite{koornwinder}, the same argument in \cite{zeilberger} can also be used when $m \ne
n$.

In this paper, I will give a new proof of the identity (\ref{eq1}) based on the Wilf-Zeilberger
method (or WZ method for short, see \cite{wilf}). This is a very succinct beautiful proof which Zeilberger
missed. I will use the following proposition (see  \cite{chen}) to prove identity
(\ref{eq1}). Given a WZ pair $(F(n,k),G(n,k))$, that is

\begin{equation}
\label{eq2}
F(n + 1,k) - F(n,k) = G(n,k + 1) - G(n,k),
\end{equation}

\noindent
then for all $m,n \in N_0 = N \cup \{0\}$, we have

\begin{equation}
\label{eq3}
\sum\limits_{k = 0}^m {F(n,k)} = \sum\limits_{j = 0}^{n - 1} {\left[ {G(j,m
+ 1) - G(j,0)} \right]} + \sum\limits_{k = 0}^m {F(0,k)}.
\end{equation}

Now set

\[
F(n,k) = {{n+k} \choose {k}}x^k(1 - x)^{n + 1}
\]

\noindent
and

\[
R(n,k) = - \frac{k}{n +
1}, G(n,k) = R(n,k)F(n,k).
\]

It is easy to verify that $(F(n,k),G(n,k))$ is a WZ pair and that the
following statements are true.

(a) For all $j \in N_0 $, $G(j,0) = 0$.

(b) For all $m \in N_0 $,$\sum\limits_{k = 0}^m {F(0,k)} = 1 - x^{m + 1}$.

By (\ref{eq3}), (a), and (b), we have

\begin{eqnarray*}
\sum\limits_{k = 0}^m {F(n,k)} &=& - x^{m + 1}\sum\limits_{j = 1}^n {{m+j} \choose {j}} (1 - x)^{j} + 1 - x^{m + 1}\\ &=& - x^{m +1}\sum\limits_{j = 0}^n {{m+j} \choose {j}}(1 - x)^{j} + 1.
\end{eqnarray*}

Finally, we have

\begin{eqnarray*}
\lefteqn{(1 - x)^{n + 1}\sum\limits_{k = 0}^m {{n+k} \choose {k}}x^{k} + x^{m + 1}\sum\limits_{k = 0}^n {{m+k} \choose {k}}(1 - x)^{k}}\\
 &=& - x^{m + 1}\sum\limits_{k = 0}^n {{m+k} \choose {k}}(1 - x)^{k} + 1 + x^{m + 1}\sum\limits_{k = 0}^n
{{m+k} \choose {k}}(1 - x)^{k}\\
& =& 1,
\end{eqnarray*}

\noindent
and the proof is complete.


\bigskip

\noindent\textit{School of Mathematical Science,
South China Normal University, Guangzhou, P.R.China\\
chenyijun73@yahoo.com.cn}

\bigskip

\end{document}